\documentclass[11pt]{article}
\newtheorem{lemma}{Lemma}[section]
\newtheorem{theorem}[lemma]{Theorem}
\newtheorem{proposition}[lemma]{Proposition}
\newtheorem{corollary}[lemma]{Corollary}
\newtheorem{definition}[lemma]{Definition}
\newtheorem{remark}[lemma]{Remark}
\newtheorem{remarks}[lemma]{Remarks}
\newtheorem{example}[lemma]{Example}
\newtheorem{examples}[lemma]{Examples}

\def\square#1#2 %maak een vierkant met dimensies #1 en inhoud #2
  {\vbox
    {\hrule
     \hbox
       {\vrule
        \vbox spread#1
          {\vfil
           \hbox spread#1{\hfil#2\hfil}
           \vfil}
        \vrule}\hrule}}
\def\sq{\square{6pt}{} }

\newenvironment{proof}{{\bf Proof:}}{\hfill $ \sq $ \vskip 4mm}

\newcommand{\thlabel}[1]{\label{th:#1}}
\newcommand{\thref}[1]{Theorem~\ref{th:#1}}
\newcommand{\selabel}[1]{\label{se:#1}}
\newcommand{\seref}[1]{Section~\ref{se:#1}}

\newcommand{\prlabel}[1]{\label{pr:#1}}

\newcommand{\colabel}[1]{\label{co:#1}}
\newcommand{\coref}[1]{Corollary~\ref{co:#1}}
\newcommand{\relabel}[1]{\label{re:#1}}
\newcommand{\reref}[1]{Remark~\ref{re:#1}}
\newcommand{\exlabel}[1]{\label{ex:#1}}
\newcommand{\exref}[1]{Example~\ref{ex:#1}}
\newcommand{\delabel}[1]{\label{de:#1}}
\newcommand{\deref}[1]{Definition~\ref{de:#1}}
\newcommand{\eqlabel}[1]{\label{eq:#1}}
\newcommand{\eqref}[1]{(\ref{eq:#1})}

\def\csh{\mbox{$\triangleright \hspace{-1ex} <$}}

\newcommand{\Hom}{{\rm Hom}}

\newcommand{\Aut}{{\rm Aut}}

\newcommand{\Char}{{\rm Char}\,}

\def\lan{\langle}
\def\ran{\rangle}

\def\pijldubbel{\lower.2ex\vbox{\hbox{$\longrightarrow$}\vspace*{-4mm}
    \hbox{$\longrightarrow$}}}
\def\doubleright#1{{\lower.2ex\vbox{\hbox{${\smash{{\cal
M}athop{\longrightarrow}
\limits^{#1}}}$}\vspace*{-4mm}\hbox{$\longrightarrow$}}}}
\def\doublerightbis#1#2{{\lower.2ex\vbox{
\hbox{${\smash{{\cal M}athop{\longrightarrow}\limits^{#1}}}$}\vspace*{-4mm}
\hbox{${\smash{{\cal M}athop{\longrightarrow}\limits_{#2}}}$}}}}
\def\doublerightleft#1#2{{\lower.2ex\vbox{
\hbox{${\smash{{\cal M}athop{\longrightarrow}\limits^{#1}}}$}\vspace*{-4mm}
\hbox{${\smash{{\cal M}athop{\longleftarrow}\limits_{#2}}}$}}}}

\def\square#1#2 %maak een vierkant met dimensies #1 en inhoud #2
{\vbox{\hrule\hbox{\vrule\vbox spread#1
{\vfil\hbox spread#1{\hfil#2\hfil}\vfil}\vrule}\hrule}}

\def\sq{\square{5pt}{} }

\def\ot{\otimes}

\def\doublerightleft#1#2{{\lower.2ex\vbox{
\hbox{${\smash{\mathop{\longrightarrow}\limits^{#1}}}$}\vspace*{-4mm}
\hbox{${\smash{\mathop{\longleftarrow}\limits_{#2}}}$}}}}

\def\ul{\underline}

\input diagrams

\begin{document}
\title{Factorisation structures of algebras and coalgebras}
\author{S. Caenepeel \\University of Brussels, VUB\\ Faculty of Applied
Sciences
\\Pleinlaan 2 \\ B-1050 Brussels, Belgium
\and Bogdan Ion\thanks{Tempus visitor at UIA.}\\
Department of Mathematics\\
 Princeton University\\
Fine Hall, Washington Road\\
Princeton, NJ 08544-1000, USA
\and G. Militaru\thanks{Research supported by the bilateral project
``Hopf algebras and co-Galois
theory" of the Flemish and Romanian governments.}
\\University of Bucharest
\\Faculty of Mathematics\\Str. Academiei 14\\RO-70109 Bucharest 1, Romania
\and Shenglin Zhu\thanks{Research supported by
the ``FWO Flanders research network WO.011.96N".}\\Institute of
Mathematics\\Fudan
University\\ Shanghai 200433, China}

\date{}
\maketitle
\begin{abstract}
\noindent
We consider the factorisation problem for bialgebras: when a bialgebra $K$
factorises as $K=HL$, where $H$ and $L$ are algebras and coalgebras
(but not necessarily bialgebras). 
Given two maps $R:\ H\ot L\to L\ot H$
and $W:\ L\ot H\to H\ot L$, we introduce a product
$L~ {}_W\hspace*{-1.5mm}\bowtie_R ~H$, and we give necessary
and sufficient conditions for $L~ {}_W\hspace*{-1.5mm}\bowtie_R ~H$ 
to be a bialgebra. It turns out that $K$ factorises as $K=HL$ if and only if
$K\cong L~ {}_W\hspace*{-1.5mm}\bowtie_R ~H$ for some maps $R$ and $W$. As examples
of this product we recover constructions introduced by Majid (\cite{Ma}) and
Radford (\cite{R0}). Also some of the pointed Hopf algebras that
were recently constructed by Beattie, D\u asc\u alescu and
Gr\"unenfelder \cite{BDG1} appear as special cases.
\end{abstract}

\section*{Introduction}\selabel{0}

The factorisation problem for a "structure" (group, algebra, coalgebra,
bialgebra) can be roughly stated as follows: in which conditions an object
$X$ can be written as a product of two subobjects $A$ and $B$ which have minimal
intersection (for example $A\cap B=\{1_X\}$ in the group case). A related 
problem is that of the construction of a new object (let us denote it by
$AB$) out of the objects $A$ and $B$. In the constructions of this type 
existing in the literature (\cite{R0}, \cite{Ta}, \cite{Ma}), the object $AB$
factorises into $A$ and $B$. This is the case - to give an example - with 
Majid's double crossed product of Hopf algebras. Moreover, whenever a Hopf 
algebra factorises in a natural way into two sub-Hopf algebras, it is 
likely to be a double crossed product (\cite[Thm. 7.2.3]{MajB}). This examples
include the quantum double of V.G. Drinfel'd, and lead also to natural 
generalisations of it on a pairing or skew-pairing of bialgeras.

The simplest example of algebra factorisation is the tensor product of
two $k$-algebras or, more general, the tensor
 product of
two algebras $A$ and $B$ in a braided monoidal category. If $A$ and $B$ are such 
algebras, 
the multiplication is then given by the
formula
\begin{equation}\eqlabel{0.1}
m_{A\# B}=(m_A\ot m_B)\circ (I_A\ot R_{A,B}\ot I_B)
\end{equation}
where $R$ is the braiding on the category. In order to apply
\eqref{0.1} to two particular algebras $A$ and $B$, we do not
need a braiding on the whole category, it suffices in fact to have
a map $R:\ B\ot A\to A\ot B$. This new algebra $A\#_R B$ will be called
a smash product, if it is associative with unit $1_A\# 1_B$.
We will work in the category of vector spaces over a field $k$, and
give necessary and sufficient conditions for $R$ to define
a smash product. The smash product can be determined completely by
a universal property, and
it will also turn out that any algebra which factorises into $A$ and $B$ is 
isomorphic with such a smash product. Therefore, we recover in this way 
several constructions that appeared earlier in the literature as special 
cases of the smash product.

The construction can be dualized, leading to the definition of the
smash coproduct of two coalgebras $C$ and $D$ (\seref{3}). The
main result of this note (\seref{4}) is the fact that we can
combine the two constructions, and this leads to the definition of
the smash biproduct of two vector space $H$ and $L$ that are
at once algebras and coalgebras (but not necessarily bialgebras).
The smash biproduct can be also characterized by what we called 
bialgebra factorisation
structures (\thref{4.3}).
Adopting this point of view for the  constructions due
to Majid \cite{Ma} and Radford \cite{R0}, we find that their constructions
are characterized by special classes of bialgebra factorisations.
We also prove that some of the pointed Hopf algebras that
Beattie, D\v asc\v alescu and Gr\"unenfelder \cite{BDG1} constructed using
iterated
Ore extensions (including classical examples like Sweedler's
four dimensional Hopf algebra) can be viewed as smash biproducts
(\thref{5.1}).

As we indicated at the beginning of this introduction, the smash biproduct
can be defined in an arbitrary monoidal category. Bernhard Drabant kindly
informed us that this general biproduct has been introduced recently
by Bespalov and Drabant in the forthcoming \cite{BD}, where it is
called a cross product bialgebra.

\section{Notations}\selabel{1}

Let $k$ be a field. For two vector spaces $V$ and $W$ and a $k$-linear
map $R:\ V\ot W\to W\ot V$ we write
$$R(v\ot w)=\sum {}^Rw\ot{}^Rv$$
for all $v\in V$, $w\in W$. Using this notation, the $k$-linear map
$$(R\ot I_V)(I_V\ot R): V\ot V\ot W\to W\ot V \ot V$$
can be denoted as follows (we write $R=r$):
$$(R\ot I_V)(I_V\ot R)(v_1\ot v_2\ot w)=
\sum {}^r(^Rw)\ot{}^rv_1\ot{}^Rv_2$$
for all $v_1,v_2\in V$, $w\in W$.

Let $A$ be a $k$-algebra. $m_A:\ A\ot A\to A$ will be the multiplication map
on $A$ and $1_A$ the unit of $A$. For a $k$-coalgebra $C$, $\Delta_C:\
C\to C\ot C$ will be the comultiplication and $\varepsilon_C:\ C\to k$
the augmentation map.

\section{The factorisation problem for algebras}\selabel{2}

Let $H,K$ and $G$ be groups. We say that $G$ factorises as $G=HK$ if $H,K$ are 
subgroups of $G$ and $H\cap K=\{1_G\}$. The problem of group factorisations was
considered before in \cite{Ta} where it led to the definition of bismash
products of Hopf algebras. One of the results regarding group factorisations
is that whenever $G$ factorises as $G=HK$, $(H,K)$ is a matched pair of groups
and $G\cong H\bowtie K$, the product associated with this pair. In order to 
prove similar results at the algebra level, we need the following definitions.
\begin{definition}
Let $A$, $B$ and $X$ be $k$-algebras with unit. We say that $X$ factorises as
$X=AB$ if there exists algebra morphisms
$$\begin{diagram}
A & \rTo^{i_A} & X & \lTo^{i_B} & B
\end{diagram}$$
such that the $k$-linear map
$$\zeta=m_X\circ(i_A\ot i_B):\ A\ot B \to X$$
is an isomorphism of vector spaces.
\end{definition}

Let $A$ and $B$ be associative $k$-algebras with unit, and
consider a $k$-linear map $R:\ B\ot A\to A\ot B$.
By definition $A\#_R B$ is equal to $A\ot B$ as a $k$-vector space
with multiplication given by the formula
\begin{equation}\eqlabel{2.1.1}
m_{A\#_RB}=(m_A\ot m_B)(I_A\ot R\ot I_B)
\end{equation}
or
\begin{equation}\eqlabel{2.1.2}
(a\#_R b)(c \#_R d)=\sum a{}^Rc\#_R {}^Rb d
\end{equation}
for all $a,c \in A$, $b, d \in B$.

\begin{definition}\delabel{2.1}
Let $A$ and $B$ be $k$-algebras with unit, and $R:\ B\ot A\to A\ot B$
a $k$-linear map. If $A\#_R B$ is an associative $k$-algebra with unit
$1_A\# 1_B$,
 we call $A\#_R B$ a smash product.
\end{definition}

\begin{remark}\relabel{2.2}\rm
Let $R$ be a braiding on the category ${\cal M}_k$ of $k$-vector spaces.
Then we have maps $R_{X,Y}:\ Y\ot X\to X\ot Y$ for all vector spaces $X$
and $Y$
and our smash product $A\#_{R_{A,B}}B$ is the usual product in the
braided category $({\cal M}_k,\ot,k,R_{X,Y})$ (see \cite{Maj}).

Our definition has a local character. In fact, given two algebras $A$ and $B$,
one can sometimes compute explicitely all the maps $R$ that make
$A\#_R B$ into a smash product (see \exref{2.11}, 3).
\end{remark}

\begin{examples}\exlabel{2.3}\label{tra}
\rm
1) Let $R=\tau_{B,A}:\ B\ot A\to A\ot B$ be the switch map. Then
$A\#_RB= A\ot B$ is the usual tensor product of $A$ and $B$.

2) Let $G$ be a group acting on the $k$-algebra $A$. This means that
we have a group homomorphism $\sigma:\ G\to \Aut_k(A)$. Writing
$\sigma(g)(a)={}^ga$, we find a $k$-linear map
$$R:kG\ot A\to A\ot kG, ~~~R(g\ot a)={}^ga\ot g$$
and $A\#_Rk[G]=A*_{\sigma}G$ is the usual skew group algebra.

3) More generaly,
let $H$ be a Hopf algebra, $A$ a left $H$-module algebra and $D$ be a 
left
$H$-comodule algebra. Let
$$
R:\ D\ot A\to A\ot D,~~~R(d\ot a)=\sum d_{<-1>}\cdot a\ot d_{<0>}.
$$
Then $A\#_R D=A\# D$ is Takeuchi's smash product \cite{Tak}.
For $D=H$, we obtain the usual smash product $A\# H$ defined in 
Sweedler's book
\cite{S}.

4) Let $(G,H)$ be a matched pair of groups and $H\bowtie G$ the product
associated to this pair (see \cite{Ta}). Write
\begin{eqnarray*}
G\times H\to H&&(g,h)\mapsto g\cdot h\\
G\times H\to G&&(g,h)\mapsto g^h
\end{eqnarray*}
for the respective group actions, and define
$$R:\ kG\ot kH\to kH\ot kG,~~~R(g\ot h)=g\cdot h\ot g^h$$
for all $g\in G$ and $h\in H$. Then $kH\#_R kG=k[H\bowtie G]$.
In \exref{2.11}, 4), we will construct an example of a smash product
$kH\#_R kG$ which is not of the form $k[H\bowtie G]$.
\end{examples}

\begin{definition}\delabel{2.4}
Let $A$ and $B$ be $k$-algebras and $R:\ B \ot A\to A\ot B$ a $k$-linear 
map.

$R$ is called left normal if
\begin{eqnarray*}
{\rm (LN)}~~~~~~&&R(b\ot 1_A)=1_A\ot b
\end{eqnarray*}
for all $b\in B$. $R$ is called right normal if
\begin{eqnarray*}
{\rm (RN)}~~~~~~&&R(1_B\ot a)=a\ot 1_B
\end{eqnarray*}
for all $a\in B$. We call $R$ normal if $R$ is left and right normal.
\end{definition}

The problem of algebra factorisations was studied before in the first part of the 
proof of Theorem 7.2.3 in \cite{MajB}. For reader's convenience we present 
in the \thref{2.5} and \thref{2.10} the detailed proofs of the results contained
there. In the following Theorem, we give necessary and sufficient conditions
for $A\#_R B$ to be a smash product.

\begin{theorem}\thlabel{sma}\thlabel{2.5}
Let $A,B$ be two algebras and let $R:B \ot A\to A\ot B$ be a $k$-linear
 map.
The following statements are equivalent
\begin{enumerate}
\item $A\#_R B$ is a smash product.
\item The following conditions hold:\\
(N) $R$ is normal;\\
(O) the following octogonal diagram is commutative
$$\begin{diagram}
B\ot A\ot A\ot B & \rTo^{I_B\ot m_A\ot I_B} & B\ot A\ot B &
                            \rTo^{R\ot I_B} & A\ot B\ot B    \\
\uTo^{I_B\ot I_A \ot R} & & & &  \dTo_{I_A\ot m_B}  \\
B\ot A\ot B\ot A &     &     &                         & A\ot B  \\
\dTo^{R\ot I_B\ot I_A} & & & &   \uTo_{m_A\ot I_B}    \\
A\ot B\ot B\ot A & \rTo^{I_A\ot m_B\ot I_A} & A\ot B\ot A &
                            \rTo^{I_A\ot R} & A\ot A\ot B
\end{diagram}$$
\item The following conditions hold:\\
(N) $R$ is normal;\\
(P) the following two pentagonal diagrams are commutative:\\
(P1)
$$\begin{diagram}
B\ot B\ot A      &\rTo^{m_B\ot I_A}&B\ot A & \rTo^R & A\ot B \\
\dTo^{I_B\ot R} &                   &       &      & \uTo_{I_A\ot m_B} \\
B\ot A \ot B      &                   & \rTo^{R\ot I_B} &    & A\ot B\ot B
\end{diagram}$$
(P2)
$$
\begin{diagram}
B\ot A\ot A      &\rTo^{I_B\ot m_A}&B\ot A & \rTo^R & A\ot B \\
\dTo^{R\ot I_A} &                   &       &      & \uTo_{m_A\ot I_B} \\
A\ot B \ot A      &                   & \rTo^{I_A\ot R} &    & A\ot A\ot B
\end{diagram}
$$
\end{enumerate}
\end{theorem}

\begin{proof}
$3)\Rightarrow 1)$ follows from \cite[Prop. 2.2 and 2.3]{vDvK}.\\
$1)\Leftrightarrow 2)$ An easy computation shows that $1_A\#1_B$ is
a right (resp. left) unit of $A\#_R B$ if and only if $R$ is right
(resp. left) normal.\\
Let us prove that the
multiplication $m_{A\#_R B}$ is associative if and only if the octogonal
diagram (O) is commutative. Using the notation introduced in \seref{1}, we
find that the commutativity of the diagram (O) is
equivalent to the following formula
\begin{equation}\eqlabel{3}\eqlabel{2.5.1}
\sum {}^r(a_2{}^Ra_3)\ot {}^rb_1{}^Rb_2=\sum{}^Ra_2{}^ra_3\ot {}^r({}^Rb_1b_2)
\end{equation}
for all $a_2,a_3\in A$, $b_1,b_2\in B$ (where $r=R$). Now for all
$a_1,a_2,a_3\in A$ and $b_1,b_2,b_3\in B$, we have that
\begin{equation}\eqlabel{2.5.2}
(a_1\#_R b_1)((a_2\#_R b_2)(a_3\#_R b_3))=
\sum a_1{}^r(a_2{}^Ra_3)\#_R{}^rb_1{}^Rb_2b_3
\end{equation}
and
\begin{equation}\eqlabel{2.5.3}
((a_1\#_R b_1)(a_2\#_R b_2))(a_3\#_R b_3)=
\sum a_1{}^Ra_2{}^ra_3\#_R{}^r({}^Rb_1b_2)b_3.
\end{equation}
and the associativity of the multiplication follows.\\
Conversely, if $m_{A\#_R B}$ is associative, then \eqref{2.5.1} follows
after we take $a_1=1_A$ and $b_3=1_B$ in (\ref{eq:2.5.2}-\ref{eq:2.5.3}).\\
$2)\Rightarrow 3)$ Suppose that (O) is commutative, or, equivalently,
\eqref{2.5.1} holds. Taking $a_2=1$ in \eqref{2.5.1}, we find, taking
the normality of $R$ into account,
\begin{equation}\eqlabel{2.5.4}
\sum {}^r({}^Ra_3)\ot {}^rb_1{}^Rb_2=
\sum {}^ra_3\ot {}^r(b_1b_2)
\end{equation}
and this is equivalent to commutativity of (P1). Taking $b_2=1$, we find
\begin{equation}\eqlabel{2.5.5}
\sum {}^r(a_2a_3)\ot {}^rb_1=\sum {}^Ra_2\ot {}^ra_3\ot {}^r({}^Rb_1)
\end{equation}
and this is equivalent to commutativity of (P1).\\
$3)\Rightarrow 2)$ Assume that (P1) and (P2) are commutative. Applying
successively \eqref{2.5.5} and \eqref{2.5.4}, we find
\begin{eqnarray*}
\sum {}^r(a_2{}^Ra_3)\ot {}^rb_1{}^Rb_2
&=& \sum {}^Ra_2{}^r({}^Ra_3)\ot {}^r({}^Rb_1){}^Rb_2\\
&=& \sum{}^Ra_2{}^ra_3\ot {}^r({}^Rb_1b_2)
\end{eqnarray*}
and this proves \eqref{2.5.1} and the commutativity of (O).
\end{proof}

\thref{2.5} leads us to the following definition.

\begin{definition}\thlabel{2.6}
Let $A $ and $B$ be $k$-algebras, and $R:\ B \ot A\to A\ot B$ a $k$-linear
map. $R$ is called left (resp. right) multiplicative if (P1)
(resp. (P2)) is commutative, or,
equivalently, if \eqref{2.5.4} (resp. \eqref{2.5.5}) holds. $R$ is
called multiplicative if $R$ is at once left and right multiplicative.
\end{definition}

\begin{remarks}\relabel{2.7}\rm
1) The conditions (\ref{eq:2.5.4}-\ref{eq:2.5.5}) for the commutativity
of the two pentagonal diagrams will turn out to be useful for computing
explicit examples. These conditions already appeared in \cite{PW}, \cite{SK},
and \thref{2.5} assures us that they can be replaced by one single
condition \eqref{2.5.1} if $R$ is normal. A map $R$ for which the diagram
(O) is commutative could be called octogonal.

2) Let $R:\ B\ot A\to A\ot B$ be normal. Then $R$ is left multiplicative
if and only if
$$
((1_A\#_R b)(1_A\#_R d))(a\#_R 1_B)=(1_A\#_R b)((1_A\#_R d)(a\#_R 1_B)),
$$
for all $a\in A$ and $b, d\in B$.
Similarly, $R$ is right multiplicative if and only if
$$
(1_A\#_R b)((a\#_R\ 1_B)(c\#_R 1_B))=((1_A\#_R b)(a\#_R 1_B))(c\#_R 1_B),
$$
for all $a,c\in A$ and $b\in B$.
\end{remarks}

\begin{examples}\exlabel{2.8}\rm
1) Let $H$ be a bialgebra and assume that $A$ and $B$ are algebras in
the category ${}_H{\cal M}$. This means that $A$ and $B$ are left
$H$-module algebras. Take $x=\sum x^1\ot x^2\in H\ot H$ and consider
the map
$$R=R_x:\ B\ot A\to A\ot B, ~~~R(b\ot a)=\sum x^2\cdot a\ot x^1\cdot b$$
If $x$ satisfies the conditions (QT1-QT4) in the definition of a
quasitriangular bialgebra (see e.g. \cite{M},\cite{Rad}), then $R_x$ is
normal and
multiplicative. Observe that we do not need that $(H,x)$ is quasitriangular,
since we do not require (QT5).

2) Now let $H$ be a bialgebra, and let $A$ and $B$ be $H$-comodule algebras,
or algebras in the category ${\cal M}^H$. For a $k$-linear map
$\sigma:\ H\ot H\to k$, we define
$R=R_{\sigma}:\ B\ot A\to A\ot B$ by
$$R(b\ot a)=\sum \sigma(a_{<1>}\ot b_{<1>})a_{<0>}\ot b_{<0>}$$
If $\sigma$ satisfies conditions (CQT1-CQT4) in the definition of a
coquasitriangular
bialgebra (cf. \cite{M}), then $R=R_{\sigma}$ is normal and multiplicative.
Observe that it is not necessary that $(H,\sigma)$ is
coquasitriangular.
\end{examples}

\begin{remarks}\relabel{2.9}\rm
1) Let $A\#_R B$ be a smash product. Then the maps
$$i_A\ :A\to A\#_R B~~{\rm and}~~i_B:\ B\to A\#_R B$$
defined by
$$i_A(a)=a\#_R 1_B~~{\rm and}~~i_B(b)=1_A\#_R b$$
are algebra maps. Moreover
$$a\#_R b=(a\#_R 1_B)(1_A\#_R b)$$
for all $a\in A$, $b\in B$.

2) If $A\#_R B$ is a smash product, then the map
$$\zeta=m_{A\#_R B}\circ (i_A\ot i_B):\ A\ot B\to A\#_R B,~~~
\zeta(a\ot b)=a\#_R b$$
is an isomorphism of vector spaces. $R$ can be reccovered from $\zeta$
by the formula
$$R=\zeta^{-1}\circ m_{A\#_R B}\circ (i_B\ot i_A)$$
\end{remarks}

This last remark leads us to the following description of the
smash product.

\begin{theorem}\thlabel{2.10}
Let $A$, $B$ and $X$ be $k$-algebras. The following conditions are equivalent.

1) There exists an algebra isomorphism $X\cong A\#_R B$, for some
$R:\ B\ot A\to A\ot B$;

2) $X$ factorises as $X=AB$.
\end{theorem}

\begin{proof}
$1)\Rightarrow 2)$ follows from \reref{2.9}.

$2)\Rightarrow 1)$
Suposse there exist algebra morphisms
$$\begin{diagram}
A & \rTo^{i_A} & X & \lTo^{i_B} & B
\end{diagram}$$
such that the $k$-linear map
$$\zeta=m_X\circ(i_A\ot i_B):\ A\ot B \to X$$
is an isomorphism of vector spaces.
Consider
$$R=\zeta^{-1}\circ m_X\circ (i_B\ot i_A):\ B\ot A\to A\ot B$$
We will prove that $R$ is normal, multiplicative and that
$\zeta:\ A\#_R B\to X$ is an algebra isomorphism.

1) $R$ is left normal. We have to show that
$R(b\ot 1_A)=1_A\ot b$, for all $b\in B$, or, equivalently,
$$(\zeta\circ R)(b\ot 1_A)=\zeta(1_A\ot b)$$
This follows easily from the following computations:
$$
(\zeta\circ R)(b\ot 1_A)=(m_X\circ(i_B\ot i_A))(b\ot 1_A)=i_B(b)i_A(1_A)=i_B(b)
$$
and
$$
\zeta(1_A\ot b)=(m_X\circ(i_A\ot i_B))(1_A\ot b)=i_A(1_A)i_B(b)=i_B(b),
$$
In a similar way, we prove that $R$ is right normal.

2) $R$ is multiplicative. To this end, it suffices to show that
\begin{equation}\eqlabel{2.10.1}
(a\#_R b)(c\#_R d)=\zeta^{-1}(\zeta(a\#_R b)\zeta(c\#_R d)),
\end{equation}
for all $a,c\in A$, $b,d\in B$. Indeed, \eqref{2.10.1} means that
the multiplication $m_{A\#_R B}$ on $A\#_R B$ can be obtained by
translating the multiplication on $X$ using $\zeta$, and this implies
that the multiplication on $A\#_R B$ is associative, and that
$\zeta$ is an algebra homomorphism. \eqref{2.10.1} is equivalent to
\begin{equation}\eqlabel{2.10.2}
\zeta((a\#_R b)(c\#_R d))=\zeta(a\#_R b)\zeta(c\#_R d)
\end{equation}
Now $\zeta(a\#_R b)\zeta(c\#_R d)=i_A(a)i_B(b)i_A(c)i_B(d)$.
Take $x= i_B(b)i_A(c)\in X$, and write
$$\zeta^{-1}(x)=\sum x_A\ot x_B\in A\ot B$$
We easily compute that
\begin{eqnarray*}
\zeta((a\#_R b)(c\#_R d))&=& \sum i_A(ax_A)i_B(x_Bd)\\
&=& i_A(a)\left( \sum i_A(x_A)i_B(x_B)\right) i_B(c)
\end{eqnarray*}
Now
\begin{eqnarray*}
\sum i_A(x_A)i_B(x_B)&=& (m_X\circ(i_A\ot i_B)\circ \zeta^{-1})\circ
(m_X\circ(i_B\ot i_A))(b\ot c)\\
&=& (m_X\circ(i_B\ot i_A))(b\ot c)\\
&=& i_B(b)i_A(c)
\end{eqnarray*}
as $m_X\circ (i_A\ot i_B)\circ \zeta=I$. This proves \eqref{2.10.2}
and completes our proof.
\end{proof}

\begin{examples}\exlabel{2.11}\label{exz}\rm
1)  Take a $k$-algebra $A$, and let
$B=k[t]$ be the polymonial ring in one variable. Consider
two $k$-linear maps $\alpha:\ A\to A$ and $\delta:\ A\to A$,
and define $R:\ B\ot A\to A\ot B$ in such a way that
$R$ is right normal and left multiplicative, and
$$R(t\ot a)= \alpha(a)\ot t +\delta(a)\ot 1$$
for all $a\in A$. This can be done in a unique way. Moreover, $R$ is
left normal if and only if $\alpha(1_A)=1_A$ and $\delta(1_A)=0$,
and $R$ is right multiplicative if and only if
$\alpha$ is an algebra morphism and $\delta$ is an $\alpha$-derivation.
If $R$ is normal and multiplicative, then $A\#_R B\cong
A[t,\alpha,\delta]$, the Ore extension associated to
$\alpha$ and $\delta$.

2) Take $a,b\in k$, and let $A=k[X]/(X^2-a),~B=k[X]/(X^2-b)$. We write
$i,j$ for the images of $X$ in respectively $A$ and $B$, and define
$R:\ B\ot A\to A\ot B$ such that $R$ is normal and
$$R(j\ot i)= -i\ot j$$
Then $A\#_R B={}^ak^b$ is nothing else then
the generalized quaternion algebra.

3) Let $C_2$ be the cyclic group of two elements. We will describe all the
smash products of the type $kC_2\#_RkC_2$. We write $A=kC_2=k[a]$ and
$B=kC_2=k[b]$ for respectively the first and second factor.
If $R:\ B\ot A\to A\ot B$ is normal, then $R$ is completely determined
by
$$R(b\ot a)= \alpha (a\ot b) + \beta (a\ot 1)+\gamma (1\ot b) +\delta (1\ot 1)$$
with $\alpha,\beta,\gamma,\delta\in k$. It is easy to check that $R$ is
multiplicative if and only if the following conditions are satisfied:
$$
2\alpha\beta=0,\quad \alpha^2+\beta^2=1,\quad
\alpha\delta+\beta\gamma+\delta=0,\quad
\alpha\gamma+ \beta\delta +\gamma=0$$
$$2\alpha\gamma=0,\quad \alpha\delta+\beta\gamma+\delta=0,\quad
\alpha^2+\gamma^2=1,\quad
\alpha\beta+ \gamma\delta+ \beta=0$$
An elementary computation shows that we have only the following possibilities
for the map $R$:
\begin{enumerate}
\item[(a)] If $\Char(k)=2$, then
   \begin{enumerate}
      \item[(i)] $R(b\ot a)= a\ot b + \delta(1\ot 1)$,
      \item[(ii)] $R(b\ot a)= (\beta+1)(a\ot b)+\beta(a\ot 1)+\beta(1\ot b)+
\beta(1\ot 1);$
   \end{enumerate}
\item[(b)] If $\Char(k)\neq 2$, then
   \begin{enumerate}
      \item[(i)] $R(b\ot a)= a\ot b,$
      \item[(ii)] $R(b\ot a)= -(a\ot b)+\delta(1\ot 1),$
      \item[(iii)] $R(b\ot a)= (a\ot 1)+(1\ot b)-(1\ot 1),$
      \item[(iv)] $R(b\ot a)= (a\ot 1)-(1\ot b)+(1\ot 1),$
      \item[(v)] $R(b\ot a)= -(a\ot 1)+(1\ot b)+(1\ot 1),$
      \item[(vi)] $R(b\ot a)= -(a\ot 1)-(1\ot b)-(1\ot 1).$
   \end{enumerate}
\end{enumerate}
We can show that the algebras given by the last four maps are isomorphic, thus
our construction gives us only the following algebras:
\begin{enumerate}
\item[(a)] If $\Char(k)=2$, then
   \begin{enumerate}
      \item[(i)] $A\#_R B=k\lan a,b|a^2=b^2=1,~ab+ba=q\ran$, where $q\in k$,
      \item[(ii)] $A\#_R B=k\lan a,b|a^2=b^2=1,~ba=(q+1)ab+qa+qb+q\ran $, where
$q\in k$;
   \end{enumerate}
\item[(b)] If $\Char(k)\not = 2$, then
   \begin{enumerate}
      \item[(i)] $A\#_R B=k\lan a,b|a^2=b^2=1,~ab=ba\ran$,
      \item[(ii)] $A\#_R B=k\lan a,b|a^2=b^2=1,~ab+ba=q\ran$, where $q\in k$,
      \item[(iii)] $A\#_R B=k\lan a,b|a^2=b^2=1,~ba=a+b-1\ran .$
   \end{enumerate}
\end{enumerate}
4) As a special case of the foregoing example, consider the normal map
$R$ defined by
$$R(b\ot a)=1_A\ot 1_B - a\ot b$$
$R$ is multiplicative, and
$$A\#_R B=k\lan a,b|a^2=b^2=1,~ab+ba=1\ran$$
is a four dimensional noncommutative $k$-algebra. Thus there exists no
group $G$ such that
$$kC_2\#_R kC_2 \cong kG, $$
and, in particular,
there exists no matched pair on $(C_2,C_2)$ such that
$kC_2\#_R kC_2 \cong k[C_2\bowtie C_2]$.\\
5) Let $A=k\lan x|x^2=0\ran \cong k[X]/(X^2)$ and $B=kC_2$ with
$C_2=\lan g\ran$. Take the unique normal map
$R:\ B\ot A\to A\ot B$ such that
$$R(g\ot x)= - x\ot g$$
Then $A\#_R B=k\lan g,x|g^2=1,~x^2=0,~gx+xg=0\ran $
is Sweedler's four dimensional Hopf algebra, considered as an
algebra.
\end{examples}

Our next aim is to show that our smash product is determined by
a universal property.

\begin{proposition}\prlabel{2.12}
Consider two $k$-algebras $A$ and $B$, and let
$R:\ B\ot A\to A\ot B$ be a normal and multiplicative map.
Given a $k$-algebra $X$, and algebra morphisms $u:\ A\to X$, $v:\ B\to X$
such that
\begin{equation}\eqlabel{2.12.1}
m_X(v\ot u)= m_X(u\ot v)R,
\end{equation}
we can find a unique algebra map $w:\ A\#_R B\to X$ such that the
following diagram commutes
$$
\begin{diagram}
  &          &A\#_R B        &          &  \\
  &\NE^{i_A} &               &\NW^{i_B} &  \\
A &          & \dTo_{w}      &          & B\\
  &\SE_{u}   &               &\SW_{v}   &  \\
  &          &  X            &          &
\end{diagram}
$$
\end{proposition}

\begin{proof}
Assume that $w$ satisfies the requirements of the Proposition. Then
$$
w(a\#_R b)=w((a\#_R 1_B)(1_A\#_R b))=(w\circ i_A)(a)(w\circ i_B)(b)=u(a)v(b),
$$
and this proves that $w$ is unique. The existence of $w$ can be proved
as follows: define $w:\ A\#_R B\to X$ by
$$w(a\#_R b)=u(a)v(b)$$
Then
$$w((a\#_R b)(c\#_R d))=\sum u(a)u({}^Rc)v({}^Rb)v(d)$$
and
$$w(a\#_R b)w(c\#_R d)=u(a)v(b)u(c)v(d)$$
and it follows from \eqref{2.12.1} that $w$ is an algebra map. The
commutativity of the diagram is obvious.
\end{proof}
%%%%%%%%%%%%%%%%%%%%%%%%%%%%%%%%%%%%%%%%%%%%%%%%%%%%%%%%%%%%%%%%%%%%%%%
%%%%%%%%%%%%%%%%%%%%%%%%%%%%%%%%%%%%%%%%%%%%%%%%%%%%%%%%%%%%%%%%%%%%%%%
%%%%%%%%%%%%%%%%%%%%%%%%%%%%%%%%%%%%%%%%%%%%%%%%%%%%%%%%%%%%%%%%%%%%%%%
\section{The coalgebra case}\selabel{3}
The results of \seref{2} can be dualized to the coalgebra case. We omit
the proofs since they are dual analogs of the corresponding proofs in
\seref{2}.

\begin{definition}
Let $C$, $D$ and $Y$ be $k$-coalgebras with counit. 
We say that $Y$ factorises as
$Y=CD$ if there exists coalgebra morphisms
$$\begin{diagram}
C & \lTo^{p_C} & Y & \rTo^{p_D} & D
\end{diagram}$$
such that the $k$-linear map
$$\eta:\ Y \to C\ot D,~~~\eta=(p_C\ot p_D)\Delta_Y$$
is an isomorphism of vector spaces.
\end{definition}

Let $C$ and $D$ be two $k$-coalgebras, and consider a $k$-linear map
$$W:\ C\ot D\to D\ot C,~~~W(c\ot d)=\sum {}^Wd\ot {}^Wc$$
Let $C\,{}_W \csh D$ be equal to $C\ot D$ as a $k$-vector space, but
with comultiplication given by
\begin{equation}\eqlabel{3.1.1}
\Delta_{C \,{}_W \csh D}= (I_C\ot W\ot I_D)(\Delta_C\ot \Delta_D)
\end{equation}
or
\begin{equation}\eqlabel{3.1.2}
\Delta_{C\,{}_W \csh D}(c \csh d)=\sum (c_{(1)} \csh {}^Wd_{(1)})
\ot ({}^Wc_{(2)} \csh d_{(2)})
\end{equation}

\begin{definition}\delabel{3.1}
With notation as above, $C \,{}_W \csh D$ is called a  smash coproduct if
the comultiplication \eqref{3.1.1} is coassociative, with counit map
$\varepsilon_{C {}_W \csh D}(c\;\csh d)=\varepsilon_C(c)\varepsilon_D(d)$.
\end{definition}

\begin{examples}\exlabel{3.2}\rm
1) Let $W=\tau_{C,D}:\ C\ot D\to D\ot C$ be the switch map.
Then $C {}_W \csh D= C\ot D$ is the usual tensor product of coalgebras.

2) Let $H$ be a bialgebra, $C$ a right $H$-module coalgebra and $D$
a right $H$-comodule coalgebra, with coaction
$\rho_D:\ D\to D\ot H,~~\rho_D=\sum d_{<0>}\ot d_{<1>}\in D\ot H$.
Let
$$W:\ C\ot D\to D\ot C,~~~W(c\ot d)=\sum d_{<0>}\ot c\cdot d_{<1>}$$
Then $C {}_W \csh D = C \csh D$ is Molnar's smash coproduct \cite{Mo}.

3) If $C$ and $D$ are finite dimensional, then $W^*$ can be viewed as
a map $W^*:\ D^*\ot C^* \to C^*\ot D^*$, after we make the identification
$(C\ot D)^*\cong C^*\ot D^*$. We have an algebra isomorphism
$$(C\, {}_W \csh D)^* \cong C^* \#_{W^*} D^*.$$
\end{examples}

\begin{definition}\delabel{3.3}
Let $C$ and $D$ be $k$-coalgebras and $W:\ C \ot D\to D\ot C$ a $k$-linear
map.\\
$W$ is called left conormal if
\begin{eqnarray*}
{\rm (LCN)}~~~~~~&&(I_D\ot \varepsilon_C)W(c\ot d)=\varepsilon_C(c)d
\end{eqnarray*}
for all $c\in C$, $d\in D$. $W$ is called right conormal if
\begin{eqnarray*}
{\rm (RCN)}~~~~~~&&(\varepsilon_D\ot I_C)W(c\ot d)=\varepsilon_D(d)c
\end{eqnarray*}
for all $c\in C$ and $d\in D$. We call $W$ conormal if $W$ is left and right
conormal.
\end{definition}

\begin{theorem}\thlabel{3.4}
Let $C$ and $D$ be $k$-coalgebras. For a $k$-linear map
$W:\ C \ot D\to D\ot C$, the following statements are equivalent.
\begin{enumerate}
\item $C {}_W \csh D$ is a smash coproduct.
\item The following conditions hold:\\
(CN) $W$ is conormal;\\
(CO) the following octogonal diagram is commutative
$$\begin{diagram}
C\ot D\ot D & \rTo^{W\ot I_D} & D\ot C\ot D &
              \rTo^{I_D\ot \Delta_C\ot I_D} & D\ot C\ot C\ot D  \\
\uTo^{I_C\ot \Delta_D} & & & &  \dTo_{I_D\ot I_C\ot W}  \\
C\ot D &     &     &                         & D\ot C\ot D\ot C  \\
\dTo^{\Delta_C\ot I_D} & & & &   \uTo_{W\ot I_D\ot I_C}    \\
C\ot C\ot D & \rTo^{I_C\ot W} & C\ot D\ot C &
              \rTo^{I_C\ot \Delta_D\ot I_C} & C\ot D\ot D\ot C
\end{diagram}$$
\item The following conditions hold\\
(CN) $W$ is conormal;\\
(CP) the following two pentagonal diagrams are commutative:\\
(CP1)
$$\begin{diagram}
C\ot D         &\rTo^{W} & D\ot C & \rTo^{\Delta_D\ot I_C} & D\ot D\ot C \\
\dTo^{I_C\ot \Delta_D} &         &         &      & \uTo_{I_D\ot W} \\
C\ot D\ot D    &       & \rTo^{W\ot I_D} &      & D\ot C\ot D
\end{diagram}$$
(CP2)
$$\begin{diagram}
C\ot D   &\rTo^{W}& D\ot C & \rTo^{I_D\ot \Delta_C} & D\ot C\ot C \\
\dTo^{\Delta_C\ot I_D} &   &       &      & \uTo_{W\ot I_C} \\
C\ot C\ot D  &                   & \rTo^{I_C\ot W} &    & C\ot D\ot C
\end{diagram}$$
\end{enumerate}
\end{theorem}

\begin{definition}\thlabel{3.5}
Let $C $ and $D$ be $k$-coalgebras, and $W:\ C \ot D\to D\ot C$ a $k$-linear
map. $W$ is called left (resp. right) comultiplicative if (CP1)
(resp. (CP2)) is commutative, or,
equivalently,
\begin{equation}\eqlabel{3.5.1}
\sum ({}^Wd)_{(1)} \ot ({}^Wd)_{(2)}\ot {}^Wc=
\sum {}^Wd_{(1)}\ot {}^wd_{(2)}\ot {}^w({}^Wc)
\end{equation}
respectively
\begin{equation}\label{3.5.2}
\sum {}^Wd\ot ({}^Wc)_{(1)} \ot ({}^Wc)_{(2)}=
\sum {}^w({}^Wd)\ot {}^wc_{(1)}\ot {}^Wc_{(2)}
\end{equation}
for all $c\in C$ and $d\in D$. $W$ is
called multiplicative if $W$ is at once left and right multiplicative.
\end{definition}

\begin{remark}\relabel{3.6}\rm
Let $W:\ C\ot D\to D\ot C$ be conormal and comultiplicative. Then the maps
\begin{eqnarray*}
p_C&:&C\,{}_W \csh D \to C,~~~p_C(c\,{}_W \csh d)=\varepsilon_D(d)c\\
p_D&:&C\,{}_W \csh D \to D,~~~p_D(c\,{}_W \csh d)=\varepsilon_C(c)d
\end{eqnarray*}
are coalgebra maps. Moreover, the $k$-linear map
$$\eta=(p_C\ot p_D)\Delta_{C\,{}_W \csh D}:\ R{}_W \csh D \to C\ot D$$
given by
$$\eta (c\,{}_W \csh d)=c\ot d$$
is bijective, and $W$ can be recovered from $\eta$ as follows:
$$W=(p_D\ot p_C)\circ\Delta_{C\,{}_W \csh D}\circ \eta^{-1}$$
\end{remark}

\begin{theorem}\thlabel{3.7}
Let $C$, $D$ and $Y$ be $k$-coalgebras. The following conditions are
equivalent.

1) There exists a coalgebra isomorphism $Y\cong C\,{}_W \csh D$, for some
$W:\ C\ot D\to D\ot C$;

2) $Y$ factorises as $Y=CD$.
\end{theorem}

Assume that $\eta$ is bijective.
From the proof of \thref{3.7}, it follows that
$$W:\ C\ot D\to D\ot C, ~~~W=(p_D\ot p_C)\circ \Delta_Y\circ \eta^{-1}$$
is conormal and commultiplicative and
$$\eta:\ Y\to C\,{}_W \csh D$$
is a isomorphism of coalgebras.\\
The smash coproduct satisfies the following universal property.

\begin{proposition}\prlabel{3.8}
Let $C $ and $D$ be $k$-coalgebras and $W:\ C\ot D\to D\ot C$ a conormal and
comultiplicative $k$-linear map.\\
Let $Y$ be a coalgebra, and $u:\ Y\to C$, $v:\ Y\to D$ coalgebra maps such that
\begin{equation}
(v\ot u)\circ \Delta_Y=W\circ (u\ot v)\circ\Delta_Y
\end{equation}
then there exists a unique coalgebra map $w:\ Y\to C{}_W \csh D$
such that $p_Cw=u$ and $p_Dw=v$:
$$\begin{diagram}
  &          &C{}_W \csh D    &          &  \\
  &\SW^{p_C} &               &\SE^{p_D} &  \\
C &          & \uTo_{w}      &          & D\\
  &\NW_{u}   &               &\NE_{v}   &  \\
  &          &  Y            &          &
\end{diagram}$$
\end{proposition}
%%%%%%%%%%%%%%%%%%%%%%%%%%%%%%%%%%%%%%%%%%%%%%%%%%%%%%%%%%%%%%%%%%%%%%%%%%%%
%%%%%%%%%%%%%%%%%%%%%%%%%%%%%%%%%%%%%%%%%%%%%%%%%%%%%%%%%%%%%%%%%%%%%%%%%%%%
%%%%%%%%%%%%%%%%%%%%%%%%%%%%%%%%%%%%%%%%%%%%%%%%%%%%%%%%%%%%%%%%%%%%%%%%%%%%

\section{Factorisations of algebras and coalgebras}\selabel{4}

Let us consider now the  problem of the factorisation of a bialgebra into 
algebras and coalgebras. We will call this problem the {\sl bialgebra
factorisation problem}.

\begin{definition}
Let $L$ and $H$ be $k$-algebras with unit which are also coalgebras
with counit and let $K$ be a bialgebra. We say that $K$ factorises as
$K=LH$ if we have maps
$$\begin{diagram}
L &\pile{\rTo^{i_L} \\ \lTo_{p_L}} & K &\pile{\lTo^{i_H} \\ \rTo_{p_H}} & H
\end{diagram}$$
such that
\begin{enumerate}
\item $i_L$ and $i_H$ are algebra maps;
\item $p_L$ and $p_H$ are coalgebra maps;
\item the $k$-linear map
$$\zeta: L\ot H \to K,~~~\zeta=m_K\circ(i_L\ot i_H)$$
is bijective and the inverse is given by
$$\zeta^{-1}: K \to L\ot H,~~~\zeta^{-1}=(p_L\ot p_H)\circ \Delta_K.$$
\end{enumerate}
\end{definition}

Let $H$ and $L$ be two algebras that are also coalgebras, and
consider two $k$-linear maps $R:\ H\ot L\to L\ot H$ and
$W:\ L\ot H\to H\ot L$. $L\, {}_W\bowtie_R \,H$ will be equal to
$L\,\#_R\, H$ as a (not necessarily associative) $k$-algebra
(see \seref{2}) and equal
to $L\, {}_W\csh \,H$ as a (not necessarily coassociative) $k$-coalgebra
(see \seref{3}).

\begin{definition}\delabel{4.1}
Let $H$, $L$, $R$ and $W$ be as above. We say that
$L~ {}_W\hspace*{-1.5mm}\bowtie_R ~H$ is a smash biproduct of $L$ and $H$
if $L\,\#_R\, H$ is a smash product, $L\, {}_W\csh \,H$ is a smash
coproduct, and $L~ {}_W\hspace*{-1.5mm}\bowtie_R ~H$ is a bialgebra.
\end{definition}

Our interest in such structures is motivated by the existing constructions 
in the Hopf algebra theory: Radford's biproducts (\cite{R0}), Takeuchi's
bismash product (\cite{Ta}) and Majid's double crossedproducts and
bicrossedproducts (\cite{Ma}). As expected, bialgebra factorisations will be 
described by the smash biproduct. Therefore, the above mentioned constructions
will become examples of biproducts.

\begin{examples}\relabel{4.2}\rm
1) Take a bialgebra $H$, and let $L$ be at once an algebra in
${}_H{\cal M}$ (an $H$-module algebra), and a coalgebra in
${}^H{\cal M}$ (an $H$-comodule coalgebra). Consider the product
$L\times H$ defined by Radford in \cite[Theorem 1]{R0}.
From \exref{2.3} 3) (with $A=L$ and $D=H$) and \exref{3.2} 2)
(with $D=L$ and $C=H$), it follows that Radford's product is
a smash biproduct.

2) In \cite{Ma}, Majid introduced the so-called bicrossproduct of two Hopf
algebras, generalizing the bismash product of Takeuchi \cite{Ta}.
We will show that Majid's bicrossproduct is an example of smash biproduct.

Let $H$ and $A$ be two Hopf algebras such that $A$ is a right $H$-module
algebra
and $H$ is a left $A$-comodule coalgebra. The structure maps are denoted
as follows
\begin{eqnarray*}
\alpha:\ A\ot H\to A&;&\alpha (a\ot h)=a\cdot h\\
\beta:\ H\to A\ot H&;&\beta(h)=\sum h_{<-1>}\ot  h_{<0>}
\end{eqnarray*}
Now consider the $k$-linear maps
\begin{eqnarray*}
R:\ A\ot H\to H\ot A&;&R(a\ot h)=\sum h_{(1)}\ot a\cdot h_{(2)}\\
W:\ H\ot A\to A\ot H&;&W(h\ot a)=\sum h_{<-1>}a\ot h_{<0>}
\end{eqnarray*}
Then $H~ _W\bowtie_R ~A$ is nothing else but Majid's
bicrossproduct $H~ ^{\beta}\bowtie_{\alpha} ~A$.

3) Majid's double crossproduct \cite[Section 3.2]{Ma} is also an
smash biproduct. To simplify notation, let us make this clear
for the non-twisted version, called bicrossed product in \cite{Ka}.

Let $(X, A)$ be a matched pair of bialgebras. This means that $X$ is a
left $A$-module coalgebra, and $A$ is a right $X$-module coalgebra
such that five additional relations hold (we refer to \cite[IX.2.2]{Ka}
for full detail). Write
\begin{eqnarray*}
\alpha:\ A\ot X\to X&;&\alpha (a\ot x)=a\cdot x\\
\beta:\ A\ot X\to A&;&\beta(a\ot x)= a^x
\end{eqnarray*}
for the structure maps, and take
\begin{eqnarray*}
R:\ A\ot X\to X\ot A&;&R(a\ot x)=\sum a_{(1)}\cdot x_{(1)} \ot
a_{(2)}^{x_{(2)}}\\
W:\ X\ot A\to A\ot X&;&W(x\ot a)=a\ot x
\end{eqnarray*}
Then $X~ _{\tau}\bowtie_R ~A$ is the double crossproduct of the matched pair
$(X,A)$ (see also \cite[Theorem IX.2.3]{Ka}).
\end{examples}

We will now give a necessary and sufficient condition for $L_W\bowtie _RH$
to be a
double product. A necessary requirement will be that
$\Delta_{L~ _W\csh ~H}$ and $\varepsilon_{L~ _W\csh ~H}$ are algebra maps,
and this will be equivalent to some compatibility relations between
$R$ and $W$. We will restrict attention to the case where $L$ and $H$
are bialgebras (see \thref{4.5}). We will first describe smash biproducts
in terms of factorisation structures.
Using \thref{2.10} and \thref{3.7}, we obtain
the following.

\begin{theorem}\thlabel{4.3}\thlabel{hex}
Let $K$ be a bialgebra and $L$, $H$ algebras that are also coalgebras.
The following statements are equivalent.

1) There exists a bialgebra isomorphism $K\cong L~ _W\bowtie_R ~H$,
for some $R:\ H\ot L\to L\ot H$ and $W:\ L\ot H\to H\ot L$;

2) There is a bialgebra factorisation of $K$ as $K=LH$.
\end{theorem}

\begin{remark}\rm\relabel{4.4}
Let us now compare the factorisation structures in Majid's constructions with
those presented here.

In \cite[Prop. 3.12]{Ma} (see also \cite[Thm. 7.2.3]{MajB}), 
it is shown that a Hopf algebra
$K$ is a double crossproduct of a matched pair of Hopf algebras $(X,A)$
if and only if  there exists a sequence
\begin{equation}\eqlabel{4.4.1}
\begin{diagram}
X &\pile{\rTo^{i_X} \\ \lTo_{p_X}} & K &\pile{\lTo^{i_A} \\ \rTo_{p_A}} & A
\end{diagram}
\end{equation}
where\\
$\bullet$ $i_X$ and $i_A$ are injective Hopf algebra maps;\\
$\bullet$ $p_X$ and $p_A$ are coalgebra maps;\\
$\bullet$ The $k$-linear map
$$\zeta:\ X\otimes A\to K, \quad \zeta:=m_K\circ(i_X\otimes i_A)$$
is bijective, and its inverse is given by the formula
$$
\zeta^{-1}:\ K\to X\otimes A, \quad \zeta^{-1}:=(p_X\otimes p_A)\circ \Delta_K
$$
Secondly, a Hopf algebra $K$ is a bicrossed product of two Hopf algebras $X$ and
$A$  if and only if there exists a sequence \eqref{4.4.1}
such that\\
$\bullet$ $i_X$ and $p_X$ are Hopf algebra maps;\\
$\bullet$ $i_A$ is an algebra map, $p_A$ is a coalgebra map;\\
$\bullet$ The $k$-linear map
$$\zeta:\ X\otimes A\to K,\quad \zeta:=m_K\circ(i_X\otimes i_A)$$ is
bijective and its inverse is given by
$$\zeta^{-1}:\ K\to X\otimes A,\quad \zeta^{-1}:=(p_X\otimes p_A)\circ
\Delta_K$$
We refer to \cite[Theorem 1.4]{Ta} for the abelian case and
to \cite[Theorem 2.3]{Ma} for the general case.\\
In a similar way, we can describe the factorisation structures
 associated to Radford's
product (see \cite[Theorem 2]{R0}). Thus, the double crossed product and
bicrossedproduct constructions are completely classifing some particular
types of bialgebra factorisation structures.
\end{remark}

We will now present neccesary and sufficient conditions for $L {}_W\bowtie _RH$
to be a smash biproduct. For technical reasons, we restricted
attention to the situation where $H$ and $L$ are bialgebras.

\begin{theorem}\thlabel{4.5}
Let $H$ and $L$ be bialgebras. For two $k$-linear maps
$R:\ H\otimes L\to L\otimes H$ and
$W:\ L\otimes H \to H\otimes L$, the following statements
are quivalent:\\
1) $L_W\bowtie _RH$ is a double product;\\
2) the following conditions hold:\\
\hspace*{14mm}(DP1)~~$R$ is normal and two-sided multiplicative;\\
\hspace*{14mm}(DP2)~~$W$ is conormal and two-sided commultiplicative;\\
\hspace*{14mm}(DP3)~~$(\varepsilon_L\otimes \varepsilon_H)R=
\varepsilon_H\otimes \varepsilon_L$;\\
\hspace*{14mm}(DP4)~~$W(1_L\otimes 1_H)=1_H\otimes 1_L$;\\
\hspace*{14mm}(DP5)~~$W(l\otimes h) = W(l\otimes 1_H)W(1_L\otimes h)$;\\
\hspace*{14mm}(DP6)~~$\sum l_{(1)}l_{(1)}^\prime\otimes^W1_H\otimes
^W[l_{(2)}l_{(2)}^\prime] =
\sum l_{(1)} \  ^Rl_{(1)}^\prime\otimes ^{^R}[^W1_H]^U1_H\otimes
^Wl_{(2)}\ ^Ul_{(2)}^\prime$;\\
\hspace*{14mm}(DP7)~~$\sum^W[h_{(1)}h^\prime_{(1)}]\otimes
^W1_L\otimes h_{(2)}h^\prime_{(2)} =
\sum ^Wh_{(1)}\ ^Uh^\prime_{(1)}\otimes ^W1_L^{^R}(^U1_L)\otimes
^Rh_{(2)}h^\prime_{(2)}$;\\
\hspace*{14mm}(DP8)~~$\sum (^Rl)_{(1)}\otimes^{^W}[(^Rh)_{(1)}]\otimes
^{^W}[(^Rl)_{(2)}]\otimes (^Rh)_{(2)} =
\sum ^Rl_{(1)}\otimes ^{^R}[^Wh_{(1)}]^U1_H\otimes
^W1_L \ ^{^r}(^Ul_{(2)})\otimes ^rh_{(2)}$,\\
for all $l,\ l'\in L,\ h,\ h'\in H$, where $r=R$ and $U=W$.
\end{theorem}

\begin{proof}
From Theorems \ref{th:2.5} and \ref{th:3.4}, it follows that
$L_W\bowtie_RH$ is at once an associative algebra with unit and a
coassociative algebra with counit if and only if (DP1) and (DP2) hold.\\
Furthermore, $\varepsilon_{L\,{}_W\csh  H}$ is an algebra map if and
only if (DP3) holds, and (DP4) is equivalent to
$$\Delta_{L_W\csh H}(1_H\bowtie 1_L) = (1_L\bowtie 1_H)\otimes (1_L\bowtie
1_H)$$
The proof will be finished if we can show that
\begin{equation}\eqlabel{4.5.1}
\Delta_{L_W\csh H}(xy)=
\Delta_{L_W\csh H}(x)\Delta_{L_W\csh H}(y)
\end{equation}
for all $x,\ y\in L_W\bowtie_RH$ if and only if (DP5-DP8) hold.
It is clear that \eqref{4.5.1} holds for all $x,\ y\in L_W\bowtie_RH$
if and only if it holds for all
$$x,y\in \{l\bowtie 1_H \mid l\in L\} \cup \{1_L\bowtie h\mid h\in H\}$$
Now \eqref{4.5.1} holds for $x=l\bowtie 1_H$ and $y= 1_L\bowtie h$
if and only if
$$\sum l_{(1)} \otimes{}^Wh_{(1)} \otimes {}^W l_{(2)} \otimes h_{(2)} =
\sum l_{(1)} \ot \ {}^W1_H\ {}^Uh_{(1)}\otimes {}^Wl_{(2)}\
{}^U1_L\otimes h_{(2)}$$
Applying $\varepsilon_L\ot I_H\ot I_L\ot \varepsilon_H$ to both sides,
we see that this condition is equivalent to (DP5).\\
\eqref{4.5.1} for all $x=l\bowtie 1_H$ and $y= l'\bowtie 1_H$ is
equivalent to (DP6),
\eqref{4.5.1} for all $x=1_L\bowtie h$ and $y= 1_L\bowtie h'$ is
equivalent to (DP7) and
\eqref{4.5.1} for all $x=1_L\bowtie h$ and $y= l\bowtie 1_H$ is
equivalent to (DP8).
\end{proof}

Let $H$ and $L$ be bialgebras, and consider a $k$-linear map
$R:\ H\otimes L\to L\otimes H$. We let $W=\tau_{L,H}:\ L\otimes H\to
H\otimes L$
be the switch map. We call
$$L{}_W\bowtie_RH=L\bowtie_RH$$
the $R$-smash product of $H$ and $L$. \thref{4.5} takes the following
more elegant form.

\begin{corollary}\colabel{4.6}
Let $H$ and $L$ be bialgebras. For a $k$-linear map $R:\ H\otimes L\to
L\otimes H$, the following statements are equivalent:

1) $L\bowtie_RH$ is an $R$-smash product;

2) $R$ is a normal, multiplicative and a coalgebra map.

Futhermore, if $H$ and $L$ have antipodes $S_H$ and $S_L$ and
$$R\tau_{L,H}(S_H\otimes S_H)R\tau_{L,H}=S_L\otimes S_H,$$
then $L\bowtie_RH$ has
an antipode given by the formula
$$S_{L\bowtie_RH}(l\bowtie h) = \sum {}^RS_L(l)\bowtie{}^RS_H(h).$$
for all $l\in L$, $h\in H$.
\end{corollary}

\begin{proof}
The switch map $W=\tau_{L,H}$ always satisfies equations (DP4-DP7),
and (DP3) and (DP8) are equivalent to $R:\ H\otimes L \to L\otimes H$
being a coalgebra map.
\end{proof}

\begin{example}\relabel{4.7}\rm
Let $H$ be a finite dimensional Hopf algebra and $L=H^{*{\bf cop}}$,
and consider the map
$$R:\\ H\otimes H^{*{\bf cop}}\to H^{*{\bf cop}}\otimes H;~~~
R(h\otimes h^*) = \sum \lan h^*, S^{-1}(h_{(3)})?h_{(1)}\ran \otimes h_{(2)}$$
In \cite{MajB}, this map $R$ is called the {\it Schr\"odinger operator}
associated
to $H$. A routine computation shows that $R$ satisfies the condition of
\coref{4.6}, and the $R$-Smash product $H^{*{\bf cop}}\bowtie_RH$ is
nothing else then the Drinfel'd Double $D(H)$ in
the sense of Radford \cite{Rad}.
\end{example}

The dual situation is also interesting:
let $H$ and $L$ be bialgebras, and take the switch map
$R=\tau_{H,L}:\ H\otimes L\to L\otimes H$. Now we call
$$L{}_W\bowtie_R H=L{}_W\bowtie H$$
the $W$-{\it smash coproduct} of $L$ and $R$, and \thref{4.5} takes
the following form.

\begin{corollary}\colabel{4.8}
Let $H$ and $L$ be bialgebras. For a $k$-linear map
$W:\ L\otimes H\to H\otimes L$, the following statements are equivalent.

1) $L{}_W\bowtie R$ is a $W$-smash coproduct;

2) $W$ is conormal, comultiplicative and an algebra map.
\end{corollary}

\begin{proof}
For $R=\tau_{H,L}$ the conditions (DP6-DP8) are equivalent to
\begin{eqnarray*}
W(ll^\prime\otimes 1_H)&=&W(l\otimes 1_H)(W(l^\prime \otimes 1_H)\\
W(1_L\otimes hh^\prime)&=&W(1_L\otimes h)W(1_L\otimes h^\prime)\\
W(l\otimes h)&=& W(1_L\otimes h)W(l\otimes 1_H)
\end{eqnarray*}
for all $l,l^\prime\in L$, $h,\ h^\prime\in H$. This three equations together
with (DP4) and (DP5) are equivalent to $W$ being an algebra map.
\end{proof}

Doi \cite{D2} and Koppinen \cite{K} introduced the category of
{\it unified Hopf modules} or {\it Doi-Hopf modules} ${\cal M}(H)_A^C$.
Some categories that are quite distinct at first sight appear as
special cases of this category: the category of classical Hopf modules
${\cal M}_H^H$ (see \cite{D2}), the category of Yetter-Drinfel'd modules
${\cal YD}^H_H$ (see \cite{CMZ1}), the category of
Long dimodules ${\cal L}^H_H$ (see \cite{Mi2}), and many others.
We will now present an alternative way to unify these categories.

\begin{definition}\delabel{4.9}
Let $H$ be a Hopf algebra and $R:\ H\otimes H\to H\otimes H$ a $k$-linear map.
A twisted $R$-Hopf module is a $k$-module $M$ that is at once a right
$H$-module and a right $H$-comodule such that the following compatibility
relation holds:
\begin{equation}\eqlabel{4.9.1}
\rho_M(m\cdot h)=\sum m_{<0>}{}^Rh_{(1)}\otimes {}^Rm_{<1>}h_{(2)}
\end{equation}
for all $m\in M$, $h\in H$.
\end{definition}

The category of twisted $R$-Hopf module and $H$-linear $H$-colinear map will
be denoted by ${\cal M}(R)^H_H$.

\begin{examples}\exlabel{4.10}\rm
1) Let $R=\tau_{H,H}$ be the switch map. Then ${\cal M}(\tau_{H,H})^H_H$
is just the category of Hopf modules, as defined in Sweedler's book \cite{S}.

2) Let $H$ be a Hopf algebra with bijective antipode and consider the map
$R:\ H\otimes H\to H\otimes H$ given by the formula
$$R(h\otimes g)=\sum g_{(2)}\otimes S^{-1}(g_{(1)})h$$
for all $h,\ g\in H$. Then the category
${\cal M}(R)^H_H$ is the category ${\cal YD}^H_H$ of Yetter-Drinfel'd modules
(see \cite{Y} and \cite{RT}).

3) Let $R:\ H\otimes H\to H\otimes H$ be given by the formula
$$R(h\otimes g)=\varepsilon(g)1_H\otimes h$$
for all $h,g\in H$. Then the category ${\cal M}(R)^H_H$ is the category
${\cal L}^H_H$ of Long $H$-dimodules defined by Long in \cite{L}.
\end{examples}

\begin{remark}\relabel{4.11}\rm
\deref{4.9} can be generalized even further. The idea is to replace the map
$R$ by a
$k$-linear map $\psi:\ C\ot A\to A\ot C$, where $C$ is a coalgebra, and $A$
is an algebra. If $\psi$ satisfies certain natural conditions, then
the triple $(A,C,\psi)$ is called an {\sl entwining structure} (see \cite{BM}),
and one can introduce the category ${\cal M}(\psi)_A^C$ of
{\sl entwining modules}. We refer to \cite{tb} for full detail.
\end{remark}
\begin{remark}{\rm
There is a close relationship between entwining structures and factorisation
structures (see \cite[Prop. 2.7]{BM}): if the coalgebra in an entwining
structure is finite dimensional, then there is a one-to-one correspondence
between entwining and factorisation structures. Therefore, the \exref{2.11} 3)
gives a complete classification of entwining structures
$(kC_2,kC_2^*,\psi)$.}
\end{remark}
%%%%%%%%%%%%%%%%%%%%%%%%%%%%%%%%%%%%%%%%%%%%%%%%%%%%%%%%%%%%%%%%%%%%%%%
%%%%%%%%%%%%%%%%%%%%%%%%%%%%%%%%%%%%%%%%%%%%%%%%%%%%%%%%%%%%%%%%%%%%%%%
%%%%%%%%%%%%%%%%%%%%%%%%%%%%%%%%%%%%%%%%%%%%%%%%%%%%%%%%%%%%%%%%%%%%%%%
\section{Examples}\selabel{5}

In \cite{BDG1}, a large class of pointed Hopf algebras is constructed,
using iterated Ore extensions. In this Section, we will see that an
important subclass of this class of Hopf algebras can be viewed as
smash biproducts of a group algebra, and a vector space that is an algebra
and a coalgebra (but not a bialgebra).

We will use the notation introduced in \cite{CD2}. Let $k$ be an algebraically
closed field of characteristic 0, $C$ a finite abelian group,
$C^*=\Hom(C, U(k))$ its character group, and $t$ a positive integer.
Assume that the following data are given:
\begin{eqnarray*}
g&=& (g_1,\ldots,g_t)\in C^t\\
g^*&=& (g^*_1,\ldots,g^*_t)\in C^{*t}\\
n&=& (n_1,\ldots, n_t)\in {\bf N}^t
\end{eqnarray*}
Assume furthermore that $\lan g^*_l,g_l\ran$ is a primitive $n_l$-th root
of unity, and that $\lan g^*_r,g_l\ran=\lan g^*_l,g_r\ran$, for all
$r\neq l$. A pointed Hopf algebra
$$K=H(C,n,g^*,g^{-1},0)$$
is then given by the following data:
\begin{equation}
x_jc= \lan g_j^*,c\ran cx_j~~~;~~~
x_jx_k= \lan g_j^*,g_k\ran x_kx_j~~~;~~~
x_j^{n_j}=0\eqlabel{5.1.1}
\end{equation}
\begin{equation}
\varepsilon(c)=1~~~;~~~\Delta(c)=c\ot c\eqlabel{5.1.2}
\end{equation}
\begin{equation}
\varepsilon(x_i)=0~~~;~~~\Delta(x_i)=x_i\ot g_i+1\ot x_i\eqlabel{5.1.3}
\end{equation}
\begin{equation}
S(c)=c^{-1}~~~;~~~S(x_i)=-x_ig_i^{-1}\eqlabel{5.1.4}
\end{equation}
for all $c\in C$ and $i\in \{1,\ldots,t\}$.
At first glance, (\ref{eq:5.1.1}-\ref{eq:5.1.4}) are not completely the
same as (3-8) in \cite{CD2}. We recover (3-8) in \cite{CD2} after we
replace $x_i$ by $g_i^{-1}x_i=y_i$.

\begin{theorem}\thlabel{5.1}
With notation as above, the Hopf algebra $K=H(C,n,g^*,g^{-1},0)$
is a smash biproduct of the group algebra $L=kC$, and an algebra
$H$ that is also a coalgebra.
\end{theorem}

\begin{proof}
We will apply \thref{4.3}. Let $H$ be the subalgebra of $K$ generated
by $x_1,\ldots,x_n$:
$$H=k\lan x_1,\ldots,x_n~|~x_j^{n_j}=0~{\rm and}~
x_jx_k=\lan g_j^*,g_k\ran x_kx_j~{\rm for~all~}j,k=1,\ldots,t\ran$$
For $m=(m_1,\ldots,m_t)\in {\bf N}^t$, we write
$$x^m=x_1^{m_1}\cdots x_t^{m_t}$$
Then
\begin{eqnarray}
\{x^m~|~0\leq m_j<n_j,~j=1,\ldots,t\}&&{\rm is~a~basis~for~}H;\eqlabel{5.1.11}\\
\{cx^m~|~0\leq m_j<n_j,~j=1,\ldots,t,~c\in C\}&&{\rm
is~a~basis~for~}H;\eqlabel{5.1.12}
\end{eqnarray}
Now $H$ and $L$ are subalgebras of $K$, and the inclusions
$$i_H:\ H\to K~~;~~i_L:\ L\to K$$
are algebra maps. Also the map
$$p_L:\ K\to L~~;~~p_L(cx^m)=\varepsilon(x^m)c$$
is a coalgebra map. Consider the maps
\begin{eqnarray*}
p_H:\ K\to H&;&p_H(cx^m)=x^m\\
\Delta_H:\ H\to H\ot H&;&\Delta_H(x^m)=(p_H\ot p_H)(\Delta_K(x^m))\\
\varepsilon_H:\ H\to k&;&\varepsilon_H(x^m)=\delta_{m0}
\end{eqnarray*}
We now claim that
\begin{equation}\eqlabel{1.5.13}
\Delta_H\circ p_H=(p_H\ot p_H)\circ\Delta_K
\end{equation}
For all $c\in C$ and $m\in {\bf N}^t$, we have
$$\Delta_H(p_H(cx^m))=\Delta_H(x^m)=(p_H\ot p_H)(\Delta_K(x^m))$$
If we write
$$\Delta_K(x^m)=\sum_{d,e\in C}\sum_{r,s\in{\bf N}^t}
\alpha_{ders} dx^r\ot ex^s$$
then
$$\Delta_K(cx^m)=\sum_{d,e\in C}\sum_{r,s\in{\bf N}^t}
\alpha_{ders} cdx^r\ot cex^s$$
and
$$(p_H\ot p_H)(\Delta_K(x^m))=(p_H\ot p_H)(\Delta_K(cx^m))=
\sum_{d,e\in C}\sum_{r,s\in{\bf N}^t}
\alpha_{ders} x^r\ot x^s$$
proving \eqref{1.5.13}. It is clear that $p_H$ is surjective, and it follows
from \eqref{1.5.13} that $\Delta_H$ is coassociative, and that
$p_H$ is a coalgebra map. Conditions 1) and 2) of \thref{4.3}
are satisfied, and condition 3) remains to be checked. Observe
that
$$\zeta(c\ot x^m)=cx^m$$
and
$$\eta(cx^m)=(p_L\ot p_H)\Delta_K(cx^m)$$
Now
\begin{equation}\eqlabel{1.5.14}
\Delta_K(cx^m)=(c\ot c)\prod_{i=1}^t (x_i\ot g_i+1\ot x_i)^{m_i}
\end{equation}
If we multiply out \eqref{1.5.14}, and apply $p_L\ot p_H$ to both sides,
then all terms are killed in the first factor by $p_L$, except
$$(c\ot c)\prod_{i=1}^t (1\ot x_i^{m_i}=c\ot cx^m$$
and it follows that
$$\eta(cx^m)=(p_L\ot p_H)(c\ot cx^m)=c\ot x^m$$
and $\eta$ is the inverse of $\zeta$, as needed.
\end{proof}

\begin{example}\exlabel{5.2}\rm
{\bf Radford's four parameter Hopf algebras}\\
Let $n,N,\nu$ be positive integers such that $n$ divides $N$ and
$1\leq \nu<N$, and let $q$ be a primitive $n$-th root of unity.
In \cite{Rdf}, Radford introduced a Hopf algebra $H_{n,q,N,\nu}$,
related to invariants for Ribbon Hopf algebras. This Hopf algebra
is a special case of the Ore extension construction of \cite{BDG1},
in fact
$$K=H_{n,q,N,\nu}=H(C_N,r,g^*,g^{-\nu},0)$$
where $g$ is a generator of the cyclic group $C_N$ of order $N$,
$t=1$, $r$ is the order of $q^{\nu}$ in $k$ and $g^*\in C^*_N$
is defined by $\lan g^*,g\ran=q$. If we take $\nu=1$ and $r=n=N$,
then we find the Taft Hopf algebra of dimension $n^2$ (see \cite{Taft}),
and Sweedler's four dimensional Hopf algebra (see \cite{S}) if
$n=2$.\\
Now $H=k[x]/(x^2)$, and the comultiplication $\Delta_K$ is given
by
\begin{equation}\eqlabel{5.2.1}
\Delta_K(g^lx^m)=
\sum_{i=0}^m {m\choose i}_{q^\nu} g^lx^i\ot g^{l+\nu i}x^{m-i}
\end{equation}
where we used the $q$-binomial coefficients and the $q$-version
of Newton's binomium (see \cite[Prop. IV.2.1]{Ka}). The comultiplication
on $H$ takes the form
\begin{equation}\eqlabel{5.2.2}
\Delta_H(x^m)=
\sum_{i=0}^m {m\choose i}_{q^\nu} x^i\ot x^{m-i}
\end{equation}
and the counit is given by $\varepsilon_H(x^m)=\delta_{m0}$.
The maps $R$ and $W$ can be described explicitely. Using \thref{2.10}
and the fact that $xg=qgx$, we find that
$R:\ H\ot L\to L\ot H$ is given by
$$R(x^m\ot g^l)=q^{lm}g^l\ot x^m$$
and, using \thref{3.7} and \eqref{5.2.1},
$$W:\ L\ot H\to H\ot L~~;~~W(g^l\ot x^m)=x^m\ot g^{l+\nu m}$$
\end{example}

\begin{example}\exlabel{5.3}\rm
The main result of \cite{CD} is that, over an algebraically closed field
of characteristic zero, there exists exactly one isomorphism class of
pointed Hopf algebras with coradical $kC_2$, represented by
$$E(n)=H(C_2,\ul{2},\ul{g}^*,\ul{g},0)$$
where $g$ and $g^*$ are the nonzero elements of $C_2$ and its dual, and
$\ul{2},\ul{g}^*,\ul{g}$ are the $t$-tuples with constant entries
$2,g^*,g$. In $E(n)$, we have the following (co)multiplication rules:
$$\begin{array}{ccc}
x_ig=-gx_i&&x_ix_j=-x_jx_i\\
\Delta(x_i)=x_i\ot g+1\ot x_i&&S(x_i)=-x_ig
\end{array}$$
Now
$$H=k\lan x_1,\ldots,x_t~|~x_i^2=0~{\rm and}~x_ix_j=-x_jx_i~{\rm for~}
i,j=1,\ldots,t\}$$
with comultiplication
$$\Delta_H(1)=1\ot 1~~~~~~~~~\Delta_H(x_i)=1\ot x_i+x_i\ot 1$$
$$\Delta(x_ix_j)=x_ix_j\ot 1+1\ot x_ix_j+x_i\ot x_j-x_jj\ot x_i$$
A little computation based on Theorems \ref{th:2.10} and \ref{th:3.7}
shows that the maps $R$ and $W$ are given by the formulas
$$\begin{array}{ccc}
R(h\ot 1)=1\ot h&&R(1\ot g)=g\ot 1\\
R(x_i\ot g)=-g\ot x_i&& R(x_ix_j\ot g)=g\ot x_ix_j\\
W(1\ot h)=h\ot 1&& W(g\ot 1)=1\ot g\\
W(g\ot x_i)=x_i\ot 1&&W(g\ot x_ix_j)=x_ix_j\ot g.
\end{array}$$
\end{example}

{\bf Acknowledgement:} We would like to thank Tomasz Brzezi\' nski for useful
comments on the first version of this paper.

\end{document}